\documentclass[12pt]{article}
\usepackage{latexsym}
\usepackage{amsfonts}

\title{On Euclidean $t$-designs}

\author{B\'ela Bajnok\\ Gettysburg College\\ bbajnok@gettysburg.edu}

\newtheorem{thm}{Theorem}
\newtheorem{defin}[thm]{Definition}

\newtheorem{cor}[thm]{Corollary}
\newtheorem{prop}[thm]{Proposition}

\newcommand{\bc}[2]{{{#1}\choose{#2}}}
\begin{document}

\renewcommand{\today}{December 28, 2004}

\maketitle

\begin{center}

\end{center}

\begin{abstract}

A Euclidean $t$-design, as introduced by Neumaier and Seidel (1988), is a finite set ${\cal X} \subset \mathbb{R}^n$ with a weight function $w: {\cal X} \rightarrow \mathbb{R}^+$ for which $$\sum_{r \in R} W_r \overline{f}_{S_{r}} = \sum_{{\bf x} \in {\cal X}} w({\bf x}) f({\bf x})$$ holds for every polynomial $f$ of total degree at most $t$; here $R$ is the set of norms of the points in ${\cal X}$,  $W_r$ is the total weight of all elements of ${\cal X}$ with norm $r$, $S_r$ is the $n$-dimensional sphere of radius $r$ centered at the origin, and $\overline{f}_{S_{r}}$ is the average of $f$ over $S_{r}$.  Neumaier and Seidel (1988), as well as Delsarte and Seidel (1989), also proved a Fisher-type inequality $|{\cal X}| \geq N(n,|R|,t)$ (assuming that the design is antipodal if $t$ is odd).  For fixed $n$ and $|R|$ we have $N(n,|R|,t)=O(t^{n-1})$.     

In Part I of this paper we provide a recursive construction for Euclidean $t$-designs in $\mathbb{R}^n$.  Namely, we show how to use certain Gauss--Jacobi quadrature formulae to ``lift'' a Euclidean $t$-design in $\mathbb{R}^{n-1}$ to a Euclidean $t$-design in $\mathbb{R}^{n}$, preserving both the norm spectrum $R$ and the weight sum $W_r$ for each $r \in R$.  For fixed $n$ and $|R|$ this construction yields a design of size $O(t^{n-1})$; however, the coefficient of $t^{n-1}$ here is significantly greater than it is in $N(n,|R|,t)$.

A Euclidean design with exactly $N(n,|R|,t)$ points is called tight; in both of the above mentioned papers it was conjectured that a tight Euclidean design with $t \geq 4$ must be a spherical design, that is, $|R|=1$ and $w$ is constant on ${\cal X}$.  Bannai and Bannai (2003) disproved this conjecture by exhibiting an example for the parameters $(n,|R|,t)=(2,2,4)$.  In Part II of this paper we construct tight Euclidean designs for $n=2$ and every $t$ and $|R|$ with $|R| \leq \frac{t+5}{4}$.  We also provide examples for tight Euclidean designs with $(n,|R|,t) \in \{(3,2,5),(3,3,7),(4,2,7)\}$.

\end{abstract}



\section{Introduction}

The concept of Euclidean designs was introduced by Neumaier and Seidel in 1988 in \cite{NeuSei:1988a} as a generalization of spherical designs, and was subsequently studied by Delsarte and Seidel in \cite{DelSei:1989a}, Seidel in \cite{Sei:1990a} and in \cite{Sei:1994a}, and just recently by Bannai and Bannai in \cite{BanBan:Prepa}.  A Euclidean design is a finite weighted set of points in the $n$-dimensional real Euclidean space $\mathbb{R}^n$ with a certain approximation property, as explained below.  First we introduce a few notations and discuss some background. 

Let $n$ be an integer and $n \geq 2$.  We denote the spaces of real polynomials, homogeneous polynomials, and homogeneous harmonic polynomials on $n$ variables by $\mathrm{Pol}( \mathbb{R}^n)$, $\mathrm{Hom}( \mathbb{R}^n)$, and $\mathrm{Harm}( \mathbb{R}^n)$, respectively.  Often we will restrict the domain of these polynomials to a subset ${\cal Y}$ of $\mathbb{R}^n$ or their degrees to a fixed integer $s$; the corresponding polynomial spaces will be denoted by $\mathrm{Pol}_s( {\cal Y})$, etc.  

The norm of a point ${\bf x}\in \mathbb{R}^n$, denoted by $|| {\bf x} ||$, is its distance from the origin; the collection of all points with given norm $r>0$ is the sphere $S_r^{n-1}$.  Let $\sigma^n$ denote the regular surface measure on $S_r^{n-1}$.  The (continuous) \emph{average value} of a polynomial $f$ on $S_r^{n-1}$ is $$\overline{f}_{S_r^{n-1}}=\frac{1}{\sigma^n(S_r^{n-1})} \int_{S_r^{n-1}} f ({\bf x}) d \sigma^n({\bf x}),$$ where 
\begin{equation}\label{GA}
\sigma^n(S_r^{n-1}) =\frac{2 r^{n-1} \pi^{\frac{n}{2}}}{\Gamma({\frac{n}{2}})}
\end{equation} is the surface area of $S_r^{n-1}$.

Let ${\cal X} \subset \mathbb{R}^n$ be a finite set.  Suppose first that every point in ${\cal X}$ has the same norm; without loss of generality assume that ${\cal X}$ is on the unit sphere, that is, ${\cal X} \subset S_1^{n-1} =: S^{n-1}$.  The concept of spherical designs, as introduced by Delsarte, Goethals, and Seidel in 1977 in \cite{DelGoeSei:1977a}, captures those sets ${\cal X}$ for which the moments of ${\cal X}$, up to a certain degree, agree with the moments of $S^{n-1}$.  Namely, we have the following definition.

\begin{defin}[Delsarte--Goethals--Seidel, \cite{DelGoeSei:1977a}]

Let ${\cal X} \subset S^{n-1}$ be a finite set and let $t$ be non-negative integer.  We say that ${\cal X}$ is a \textup{\textbf{spherical $t$-design}} if 
$$ |{\cal X}| \overline{f}_{S^{n-1}} = \sum_{{\bf x} \in {\cal X}}  f({\bf x})$$ holds for every $f \in \mathrm{Pol}_t( \mathbb{R}^n)$.

\end{defin}

Spherical designs enjoy a vast and rapidly growing literature, and have been studied from a variety of perspectives, including algebra, combinatorics, functional analysis, geometry, number theory, numerical analysis, and statistics.  For general references, please see \cite{Ban:1988a}, \cite{ConSlo:1999a}, \cite{DelPac:Prepa}, \cite{DelGoeSei:1977a}, \cite{EriZin:2001a}, \cite{God:1993a}, \cite{GoeSei:1979a}, \cite{GoeSei:1981a}, \cite{Sei:1996a}, and \cite{SeyZas:1984a}.

Euclidean designs generalize spherical designs in two aspects: we do not assume that all points in the design have the same norm, and we allow the points to have different weights.  

For a finite set ${\cal X} \subset \mathbb{R}^n \setminus \{ {\bf 0} \}$, we call $R=\{ || {\bf x} ||  \mbox{   } | \mbox{   } {\bf x} \in {\cal X} \}$ the \emph{norm spectrum} of ${\cal X}$.  (In this paper, for convenience, we exclude the possibility of ${\bf 0} \in {\cal X}$; see \cite{BanBan:Prepa} for a discussion.)  We can partition ${\cal X}$ into \emph{layers} ${\cal X} = \cup_{r \in R} {\cal X}_{r}$ where ${\cal X}_r={\cal X} \cap S_{r}^{n-1}$.   A weight function on ${\cal X}$ is a function $w: {\cal X} \rightarrow \mathbb{R}^+$; the \emph{weight distribution of $w$ on $R$} is the function $W: R \rightarrow \mathbb{R}^+$ given by $W_r=\sum_{{\bf x} \in {\cal X}_r} w({\bf x})$.   (Throughout this paper, we only consider positive weights.)  We are now ready to state our definition of Euclidean $t$-designs.

\begin{defin}[Neumaier--Seidel, \cite{NeuSei:1988a}]

Let ${\cal X} \subset \mathbb{R}^n \setminus \{ {\bf 0} \}$ be a finite set with norm spectrum $R$, and suppose that a weight function $w$ is given on ${\cal X}$ which has weight distribution $W$ on $R$.  Let $t$ be non-negative integer.  We say that $({\cal X},w)$ is a \textup{\textbf{Euclidean $t$-design}} if 
$$\sum_{r \in R} W_r \overline{f}_{S_{r}^{n-1}} = \sum_{{\bf x} \in {\cal X}} w({\bf x}) f({\bf x})$$ holds for every $f \in \mathrm{Pol}_t( \mathbb{R}^n)$.

\end{defin}

In Part I of this paper we provide a recursive construction for Euclidean $t$-designs in $\mathbb{R}^n$.  Namely, we show how to use certain Gauss--Jacobi quadrature formulae to ``lift'' a Euclidean $t$-design $({\cal A}^{n-1},w^{n-1})$ in $\mathbb{R}^{n-1}$ to a Euclidean $t$-design $({\cal A}^{n},w^{n})$ in $\mathbb{R}^{n}$.  Our recursion will preserve both the norm spectrum $R$ and the weight distribution $W$ on $R$; that is, we will have $$R=\{ || {\bf x} ||  \mbox{   } | \mbox{   } {\bf x} \in {\cal A}^{n-1} \}=\{ || {\bf x} ||  \mbox{   } | \mbox{   } {\bf x} \in {\cal A}^{n} \}$$ and $$W_r=\sum_{{\bf x} \in {\cal A}^{n-1}_r} w({\bf x})=\sum_{{\bf x} \in {\cal A}^{n}_r} w({\bf x})$$ for each $r \in R$.  Since it is not difficult to exhibit Euclidean $t$-designs in the plane (see, for example, Part II of this paper), our recursion yields Euclidean $t$-designs in any dimension.  We can describe this construction more precisely, as follows.

First, let us define an analogue of spherical designs for a real interval.  Let us set $I=[-1,1]$ and suppose that $\mu: I \rightarrow \mathbb{R}$ is Lebesgue-integrable function with positive integral on every non-degenerate subinterval of $I$.  For a polynomial $f \in \mathrm{Pol}( \mathbb{R})$ we define the \emph{$\mu$-average of $f$} to be
$$ \overline{f}_{\mu} = \frac{1}{\mu(I)} \int_{I} \mu(x) f (x) d x$$ where $$\mu(I)=\int_I \mu(x) dx.$$

\begin{defin}

Let ${\cal C} \subset I$ be a finite set and let $\gamma :{\cal C} \rightarrow \mathbb{R}^+$ be a (weight) function.  Let $t$ be non-negative integer and let $\mu: I \rightarrow \mathbb{R}$ be a Lebesgue-integrable function with positive integral on every non-degenerate subinterval of $I$.  We say that $({\cal C},\gamma)$ is an \textup{\textbf{interval $t$-design}} on $I$ with respect to the function $\mu$ if 
$$ \overline{f}_{\mu} = \sum_{c \in {\cal C}}  \gamma(c) f(c)$$ holds for every $f \in \mathrm{Pol}_t( \mathbb{R})$.   

\end{defin}

We will discuss interval $t$-designs in more detail in Section 2.  In particular, we will focus on the case of the ultraspherical weight function, in which case the nodes are roots of certain Gegenbauer polynomials, and the weights are also expressed in terms of these roots.  Furthermore, our interval design will be \emph{interior}, that is, we will have ${\cal C} \subset (-1,1)$; and \emph{antipodal}, that is, for every $c \in {\cal C}$, we will have $-c \in {\cal C}$ and $\gamma(c)=\gamma(-c)$.  

Next, we introduce some notations for ``lifting'' points from $\mathbb{R}^{n-1}$ to $\mathbb{R}^{n}$.  For a point ${\bf y}=(y_1,\dots,y_{n-1}) \in \mathbb{R}^{n-1}$ and a number $c \in I$, let us define ${\bf y} \star c \in \mathbb{R}^{n}$ by $${\bf y} \star c :=(\sqrt{1-c^2}y_1, \dots, \sqrt{1-c^2}y_{n-1}, c||{\bf y} ||);$$  for ${\cal Y} \subseteq \mathbb{R}^{n-1}$ and ${\cal C} \subseteq I$ we set $${\cal Y} \star {\cal C} = \{ {\bf y} \star c \mbox{  } | \mbox{  } {\bf y} \in {\cal Y}, c \in {\cal C} \}.$$  Note that we always have $||{\bf y}\star c ||=||{\bf y} || ,$ and therefore $$S_r^{n-2} \star I=S_r^{n-1}.$$  

For functions $f: \mathbb{R}^{n-1} \rightarrow \mathbb{R}$ and $g: I \rightarrow \mathbb{R}$, we define $f \star g: \mathbb{R}^{n-1} \star I \rightarrow \mathbb{R}$ by $$(f \star g)({\bf y} \star c)=f({\bf y}) g(c).$$

In Section 3 we will prove the following.

\begin{thm}\label{FC}

Suppose that $({\cal C}^{n}, \gamma^{n})$ is an interior and antipodal interval $t$-design on $I$ with respect to the (ultraspherical) weight function $$\mu(x)=(1-x^2)^{\frac{n-3}{2}}.$$  Suppose that $({\cal A}^{n-1}, w^{n-1})$ is any Euclidean $t$-design in $\mathbb{R}^{n-1}$, and define $$({\cal A}^{n}, w^{n})=({\cal A}^{n-1} \star {\cal C}^{n}, w^{n-1} \star \gamma^{n}),$$ as above.  Then $({\cal A}^{n-1}, w^{n-1})$ is a Euclidean $t$-design in $\mathbb{R}^{n}$.  Furthermore, if $({\cal A}^{n-1}, w^{n-1})$ has norm spectrum $R$ and weight distribution $W$ on $R$, then $({\cal A}^{n}, w^{n})$ will also have norm spectrum $R$ and weight distribution $W$ on $R$.

\end{thm}

We note that our recursive construction for Euclidean $t$-designs is ``semi-explicit'' in the sense that the co\"{o}rdinates and the weights of the points in the design are given in terms of the roots of certain Gegenbauer polynomials; these roots can be computed explicitly, however, whenever $t \leq 17$ (see Section 2).

In Part II of the paper we investigate Euclidean designs of minimum size.  Before describing our results here, let us summarize what is known about spherical designs of minimum size.  

For a non-negative integer $s$, let 
$$d(s):=\dim \mathrm{Hom}_s( \mathbb{R}^n) = \bc{s+n-1}{n-1}.$$
We have the following \emph{Fisher-type inequality}.

\begin{thm}[Delsarte--Goethals--Seidel, \cite{DelGoeSei:1977a}]\label{FA}

Let $$N_1=d \left( \left\lfloor \frac{t}{2} \right\rfloor \right) + d \left( \left\lfloor \frac{t-1}{2} \right\rfloor \right).$$
Then for every spherical $t$-design ${\cal X} \subset S^{n-1}$ we have $|{\cal X}| \geq N_1$.    

\end{thm} 

\begin{defin}

Spherical $t$-designs of minimum size $N_1$ are called \textup{\textbf{tight}}.

\end{defin} 

Examples for tight spherical designs in the plane ($n=2$) are provided by the vertices of regular $(t+1)$-gons.  For $n \geq 3$ we see that two antipodal points, the vertices of the regular simplex, and the vertices of the generalized regular octahedron form tight $t$-designs for $t=1$, $t=2$, and $t=3$, respectively.  For $n \geq 3$ and $t \geq 4$ examples are far and between.  In particular, Bannai and Damerell proved in \cite{BanDam:1979a} and \cite{BanDam:1980a} that for $n \geq 3$ we must have $t \in \{1,2,3,4,5,7,11\}$; furthermore, for $t=11$, the only tight design, up to isometry, is provided by the Leech lattice ($n=24$ and $|{\cal X}|= 196,560$).  The possibly finite number of cases for $t=4, 5$ and 7 have not quite been classified (see the upcoming paper \cite{BanMunVen:Prepa} of Bannai, Munemasa, and Venkov for recent progress).

The corresponding Fisher-type inequality for Euclidean designs provides the minimum size for the case of an even $t$; for odd $t$ we only have a lower bound if the design is \emph{antipodal}, that is, for every $x \in {\cal X}$, $ -x \in {\cal X}$ and $w(-x)=w(x)$.

\begin{thm}[Delsarte--Seidel, \cite{DelSei:1989a}]\label{FB}  For a positive integer $k$, let $$N_k=d \left( \left\lfloor \frac{t}{2} \right\rfloor +2-2k \right) + d \left( \left\lfloor \frac{t-1}{2} \right\rfloor +2-2k \right).$$ Let $({\cal X},w)$ be a Euclidean $t$-design on $p$ layers in $\mathbb{R}^n$; if $t$ is odd, assume further that the design is antipodal.
Then we have $$|{\cal X}| \geq N(n,p,t):=\sum_{k=1}^{p} N_k.$$ 

\begin{defin}\label{JJ}

We say that a Euclidean $t$-design on $p$ concentric spheres is \textup{\textbf{tight}} if it has size $N(n,p,t)$.

\end{defin}

\end{thm}  

\emph{Remarks.}  1.  Note that $N(n,1,t)=N_1$ above agrees with $N_1$ in Theorem \ref{FA}.  Delsarte, Goethals, and Seidel proved in \cite{DelGoeSei:1977a} that, for an odd value of $t$, a tight spherical $t$-design is necessarily antipodal (it is not known whether this holds for Euclidean $t$-designs as well).  Therefore Theorem \ref{FA} can be regarded as a special case of Theorem \ref{FB}.

2.  There is another remark we make about the case $p=1$: According to Theorems \ref{FA} and \ref{FB}, a \emph{weighted spherical design} (a set where the points lie on the same sphere but are allowed to have different weights) cannot have fewer points than a tight spherical design has (with constant weight).  In fact, more is true: a result of Bannai and Bannai \cite{BanBan:Prepa} implies that, at least when $t$ is even, a tight Euclidean design with $p=1$ \emph{must} have constant weights.

3.  We also point out an observation when $p$ is relatively large.  Note that $N_k=0$ if $\left\lfloor \frac{t}{2} \right\rfloor +2-2k <0$; in particular, the minimum size $N(n,p,t)$ remains constant once $p \geq \left\lfloor \frac{t}{4} \right\rfloor +1 $.  Accordingly, Bannai and Bannai \cite{BanBan:Prepa} separate Definition \ref{JJ} into the cases $p \leq \left\lfloor \frac{t}{4} \right\rfloor$ and $p \geq \left\lfloor \frac{t}{4} \right\rfloor +1 $ by using the terms ``tight $t$-design on $p$ concentric spheres'' and ``Euclidean tight $t$-design'', respectively (\cite{BanBan:Prepa} only deals with the case when $t$ is even).

4.  The following bounds might be useful.
$$\bc{\left\lfloor \frac{t}{2} \right\rfloor +n}{n}-\bc{ \left\lfloor \frac{t}{2} \right\rfloor -2p +n}{n} \leq N(n,p,t) \leq \bc{ \left\lfloor \frac{t+1}{2} \right\rfloor +n}{n}-\bc{\left\lfloor \frac{t+1}{2} \right\rfloor -2p+n}{n};$$ in fact, when $t$ is even, we have the closed form $$ N(n,p,t) =\bc{\frac{t}{2}+n}{n}-\bc{ \frac{t}{2} -2p+n}{n}.$$
Furthermore, if $t$ is even and $p \geq \left\lfloor \frac{t}{4} \right\rfloor +1$, then $$N(n,p,t)=\bc{ \frac{t}{2}+n}{n},$$ a quantity which coincides with the maximum size of a $\frac{t}{2}$-\emph{distance set} in $\mathbb{R}^n$ (see \cite{BanBanSta:1983a}), and therefore a Euclidean $t$-design on $p \geq \left\lfloor \frac{t}{4} \right\rfloor +1$ spheres cannot be an $s$-distance set for $s<\frac{t}{2}$. 

In Section 5 of this paper we explicitly construct tight Euclidean designs in the plane for every $t$ and every $p \leq \left\lfloor \frac{t+5}{4}\right\rfloor$.  Our construction is a generalization of the example for $p=2$ and $t=4$ which already appeared in \cite{BanBan:Prepa}.  Namely, we will show how a union of concentric regular polygons, with appropriate weights and rotations about the origin, forms a tight Euclidean design.

\begin{thm}\label{FD}

Let $1 \leq p \leq \left\lfloor \frac{t+5}{4}\right\rfloor$, and set $m=t+3-2p$.  Let $R=\{r_1,\dots,r_p\}$ be an arbitrary set of $p$ distinct positive real numbers; without loss of generality assume that $r_1<r_2< \cdots < r_p$.  For integers $k$ and $j$, define the point $$b_{k,j}= \left(r_k \cos \left( \frac{2j+k}{m} \pi \right), r_k \sin \left( \frac{2j+k}{m} \pi \right) \right),$$ and let $${\cal B} = \{ b_{k,j} \mbox{  } | \mbox{  } 1 \leq j \leq m , 1 \leq k \leq p \}.$$   Furthermore, define the weight function $w : {\cal B} \rightarrow \mathbb{R}^+$ with
$$w(b_{k,j})=\left\{ \begin{array}{lll}
\frac{1}{r_1^m} & \mbox{ if } & k=1 \\
(-1)^k \frac{1}{r_k^m} \prod_{2 \leq l \leq p, l \not = k} \frac{r_1^2-r_l^2}{ r_k^2-r_l^2} & \mbox{ if } & 2 \leq k \leq p
\end{array} \right.$$
Then $({\cal B}, w)$ is a tight Euclidean $t$-design in the plane with norm spectrum $R$.  (Note that $w(b_{k,j})$ is always positive and that ${\cal B} $ is antipodal when $t$ is odd.)

\end{thm}

Combining Theorems \ref{FC} and \ref{FD} (also applying Theorem \ref{FE} of Section 2), we get Euclidean designs in any dimension.

\begin{cor}\label{FF} If $1 \leq p \leq \frac{t+5}{4}$, then there is a Euclidean $t$-design in $\mathbb{R}^n$ with norm spectrum of size $p$, consisting of
\begin{equation}\label{A}
p \cdot (t+3-2p) \cdot \left( \left\lfloor \frac{t}{2} \right\rfloor +1 \right)^{n-2}
\end{equation} 
points.
\end{cor}

\emph{Remark.}  Note that, for fixed $n$ and $p$ and when $t$ approaches infinity, (\ref{A}) is of order $O(t^{n-1})$; by Remark 4 after Theorem \ref{FB}, $N(n,p,t)$ is also of order $O(t^{n-1})$.  We note, however, that (\ref{A}) is still significantly larger than $N(n,p,t)$.  Nevertheless, this is drastically better than the best analogous result known for spherical designs; while it is conjectured that spherical $t$-designs in $n$ dimensions exist of size $O(t^{n-1})$ as well, this has only been proved for $O(t^{\frac{n^2-n}{2}})$; see \cite{KorMey:1994a}, \cite{Kui:1995a}, \cite{Baj:1991d}, and \cite{Baj:1992a}.   

Let us now turn to the case of tight Euclidean designs in higher dimensions ($n \geq 3$).  Here we restrict our search to \emph{fully symmetric designs}, that is, those which remain fixed by permutations of the co\"ordinates and by reflections with respect to co\"ordinate hyperplanes.  In particular, we consider only antipodal sets, and we may assume that $t$ is odd.  

For $t=1$ a pair of two antipodal points, with equal weights, forms a tight 1-design in $\mathbb{R}^n$; for $t=3$, the $n$ pairs $$\{ (\pm r_1, 0, \dots, 0), (0, \pm r_2, 0, \dots, 0), \dots, (0, \dots, 0, \pm r_n)\},$$ with weights inversely proportional to their squared norms, form tight 3-designs in $\mathbb{R}^n$.  

Here we provide examples of tight Euclidean designs for $n=3$ with $t=5$ and $t=7$, and for $n=4$ with $t=7$.  In particular, we will prove that 

\begin{enumerate}
\item
in $\mathbb{R}^3$, the union of an octahedron and a cube, with appropriate weights, forms a tight 5-design;
\item
in $\mathbb{R}^3$, the union of an octahedron, a cuboctahedron, and a cube, with appropriate weights, forms a tight 7-design; and
\item
in $\mathbb{R}^4$, the points of minimum non-zero norm in the lattice $D_4$ together with the points of minimum non-zero norm in the dual lattice $D_4^*$, with appropriate weights, form a tight 7-design.
\end{enumerate}

It might be useful to summarize these examples more precisely, as follows.
 
\begin{prop}\label{FG}

For $1 \leq k \leq n$, define $${\cal I}^n_k =\{ {\bf x} \in \{-1,0,1\}^n \mbox{   } | \mbox{  } ||{\bf x} ||^2=k \}.$$  For a given (index) set $J \subseteq \{1,2, \dots, n\}$ and functions $r: J \rightarrow \mathbb{R}^+$ and $w: J \rightarrow \mathbb{R}^+$, we let  $\Psi=(J, r, w)$; we then consider the set $${\cal X}={\cal X}(\Psi)=\cup_{k \in J} \frac{r(k)}{\sqrt{k}} {\cal I}^n_k$$ with weight function $w=w(\Psi): {\cal X} \rightarrow \mathbb{R}^+$ given by $$w({\bf x})=w(k) \mbox{   if   } {\bf x} \in \frac{r(k)}{\sqrt{k}}{\cal I}^n_k.$$

Let $\lambda \in \mathbb{R}^+$, $\lambda \not = 1$.

\begin{enumerate} 

\item For $n=3$, $p=2$, and $t=5$, choose $\Psi$ according to the following table.
$$\begin{array}{ccc}  
k \in J & 1 & 3  \\ 
r & 1 & \lambda \\ 
w & 1 & \frac{9}{8 \lambda^4} \\  
\end{array}$$ 
Then $({\cal X},w)$ is a tight Euclidean 5-design of size 14 in $\mathbb{R}^3$.

\item For $n=3$, $p=3$, and $t=7$, choose $\Psi$ according to the following table.
$$\begin{array}{cccc}  
k \in J & 1 & 2 & 3  \\  
r & 1 & \sqrt{\frac{2\lambda+3}{5\lambda}} & \sqrt{\frac{2\lambda+3}{5}} \\ 
w & 1 & \frac{100\lambda^3}{(2\lambda+3)^3} & \frac{675}{8(2\lambda+3)^3} \\  
\end{array}$$ 
Then $({\cal X},w)$ is a tight Euclidean 7-design of size 26 in $\mathbb{R}^3$.

\item For $n=4$, $p=2$, and $t=7$, choose $\Psi$ according to the following table.
$$\begin{array}{cccc}  
k \in J & 1 & 2 & 4  \\  
r & 1 & \lambda & 1 \\ 
w & 1 & \frac{1}{\lambda^6} & 1 \\  
\end{array}$$ 
Then $({\cal X},w)$ is a tight Euclidean 7-design of size 48 in $\mathbb{R}^4$.

\end{enumerate}

\end{prop}

We point out that both Theorem \ref{FD} and Proposition \ref{FG} disprove a conjecture of Neumaier and Seidel in \cite{NeuSei:1988a} (see also \cite{DelSei:1989a}) that there are no tight Euclidean $t$-designs with $p \geq 2$ and $t \geq 4$.  It seems to be an interesting problem to classify all other tight Euclidean designs.

\vspace*{.2in}

\begin{center}

{\Large{{\bf Part I: A Recursive Construction of Euclidean Designs}}}

\end{center}

\section{Interval designs with respect to the ultraspherical weight function}

Let us here review some information on Gauss--Jacobi quadrature formulae.

Recall that, for fixed $\alpha, \beta > -1$, the Jacobi polynomial $P_s^{(\alpha,\beta)}: I=[-1,1] \rightarrow \mathbb{R}$  can be defined as $$P_s^{(\alpha,\beta)}(x)=\frac{(-1)^s}{2^ss!(1-x)^{\alpha}(1+x)^{\beta}}\frac{d^s}{dx^s}[(1-x)^{\alpha+s}(1+x)^{\beta+s}].$$ 
The Jacobi polynomial $P_s^{(\alpha,\beta)}$ has degree $s$ and has $s$ distinct roots in the interval $(-1,1)$; these roots, with appropriate weights (sometimes called \emph{Christoffel numbers}) form a \emph{quadrature formula} of strength at least $2s-1$ with respect to the weight function $$\mu^{(\alpha,\beta)}(x)=(1-x)^{\alpha}(1+x)^{\beta}.$$ This is to say that, if ${\cal C}_s^{(\alpha,\beta)}$ denotes the set of roots of $P_s^{(\alpha,\beta)}$ and $\gamma_s^{(\alpha,\beta)}(c)$ is the weight of $c \in {\cal C}_s^{(\alpha,\beta)}$ defined as $$\gamma_s^{(\alpha,\beta)}(c)=\frac{2^{\alpha+\beta+1} \Gamma(\alpha+s+1) \Gamma(\beta+s+1)}{s! \Gamma(\alpha+\beta+s+1)} \cdot \frac{1}{(1-c^2)[(P^{(\alpha,\beta)}_s(c))']^2},$$ then 
$$\int_{I} f(x) \mu^{(\alpha,\beta)}(x) dx = \sum_{c \in {\cal C}_s^{(\alpha,\beta)}} \gamma_s^{(\alpha,\beta)}(c) f(c)$$ holds for every $f \in \mathrm{Pol}_{2s-1}(\mathbb{R})$.  

We here will need only the (ultraspherical) case when, for our fixed $n$, $\alpha=\beta=\frac{n-3}{2}$; we let $$P_s^{\frac{n-3}{2}}(x)= P_s^{(\frac{n-3}{2},\frac{n-3}{2})}(x),$$ $${\cal C}_s^n={\cal C}_s^{(\frac{n-3}{2},\frac{n-3}{2})},$$ and $$\mu(x)=(1-x^2)^{\frac{n-3}{2}}.$$  Then $P_s^{\frac{n-3}{2}}$ (called a \emph{Gegenbauer polynomial}) is an even function when $s$ is even and an odd function when $s$ is odd, and therefore its $s$ distinct roots form an antipodal set in $(-1,1)$.   As a consequence, the roots of $P_s^{\frac{n-3}{2}}$ can be, at least in theory, found explicitly for every $s \leq 9$.  Furthermore, it might be worth to point out that when $n=2$ or $n=4$, the corresponding Jacobi polynomials are the Chebyshev polynomials of the first kind and second kind, respectively; in these cases we can express their roots as $${\cal C}_s^2=\left\{ \cos \left(  \frac{2i-1}{2s} \pi \right) \mbox{  } | \mbox{  } 1 \leq i \leq s \right\}$$ and $${\cal C}_s^4=\left\{\cos \left(  \frac{i}{s+1} \pi \right) \mbox{  } | \mbox{  } 1 \leq i \leq s \right\}.$$ 

We summarize these findings using the terminology of interval designs.

\begin{thm}\label{FE}

Let $t$ be a non-negative integer, and let $s=\lfloor \frac{t}{2} \rfloor +1$.  Then there is an interior and antipodal interval $t$-design $({\cal C}_s^n, \gamma_s^n)$ on $I$ with respect to the weight function $\mu(x)=(1-x^2)^{\frac{n-3}{2}}.$  Furthermore, $|{\cal C}_s^n|=s$, and these points and their weights can be computed exactly if $n \in \{2,4\}$ (and $t$ arbitrary) or when $t \leq 17$ (and $n$ arbitrary). 

\end{thm}

For more details on these and other facts about the Gauss--Jacobi quadratures, see e.g. \cite{Hil:1974a} or \cite{Kry:1962a}.

\section{The recursive construction}

In this section we prove Theorem \ref{FC}.  

As before, let $I=[-1,1]$ and set $$\mu(x)=(1-x^2)^{\frac{n-3}{2}}.$$  Then we have 
\begin{equation}\label{GB}
\mu(I)=\int_I (1-x^2)^{\frac{n-3}{2}} dx=\frac{\sqrt{\pi} \cdot \Gamma(\frac{n-1}{2})}{\Gamma(\frac{n}{2})}.
\end{equation}

Let $({\cal C}^n, \gamma^n)$ be an interior and antipodal $t$-design in $I$ with respect to $\mu$.  (According to Theorem \ref{FE}, we may have $|{\cal C}^n|=s=\left \lfloor \frac{t}{2} \right \rfloor +1$.)  With the notations $${\cal C}^{n} = \{ c_i^{n} \mbox{   } | \mbox{  } 1 \leq i \leq |{\cal C}^n|\}$$ and $$\gamma^{n}(c_i^{n})=\gamma^{n}_i,$$ the quadrature formula 
\begin{equation}\label{AE}
\frac{1}{\mu(I)} \int_I f(x) (1-x^2)^{\frac{n-3}{2}}  d x = \sum_{i} \gamma^{n}_{i} f(c_i^n)
\end{equation} 
is then exact for every $f \in \mathrm{Pol}_t( \mathbb{R})$.  In particular, note that $$\sum_{i} \gamma^{n}_{i}=1.$$

Assume next that $({\cal A}^{n-1}, w^{n-1})$ is a Euclidean $t$-design in $\mathbb{R}^{n-1}$ with norm spectrum $$R=\{r_k \mbox{  } | \mbox{  } 1 \leq k \leq p \}.$$  By introducing notations $${\cal A}^{n-1}= \cup_{k=1}^p {\cal A}_k^{n-1},$$
$${\cal A}_k^{n-1} = \{ {\bf a_{k,j}^{n-1}}=(a_{k,j,1}^{n-1},a_{k,j,2}^{n-1}, \dots, a_{k,j,n-1}^{n-1}) \}_j \subset S_{r_k}^{n-2} \subset \mathbb{R}^{n-1},$$ $$w^{n-1}({\bf a_{k,j}^{n-1}})=w^{n-1}_{k,j},$$ and $$W_k^{n-1}=\sum_j w^{n-1}_{k,j},$$ we can write this as the identity 
\begin{equation}\label{AD}
\sum_k \frac{W_k^{n-1}}{\sigma^{n-1}(S_{r_k}^{n-2})} \int_{S_{r_k}^{n-2}} f({\bf x}) d \sigma^{n-1}({\bf x}) = \sum_{k,j} w^{n-1}_{k,j} f({\bf a_{k,j}^{n-1}})
\end{equation} 
which then holds for every $f \in \mathrm{Pol}_t(\mathbb{R}^{n-1})$.

We then define $${\bf a_{k,j,i}^{n}}=(a_{k,j,i,1}^{n},a_{k,j,i,2}^{n}, \dots, a_{k,j,i,n}^{n}) ={\bf a_{k,j}^{n-1}} \star c_i^n \in \mathbb{R}^n;$$ recall that this means $$a_{k,j,i,l}^{n}=\left\{ \begin{array}{lll}
\sqrt{1-(c_i^n)^2} a_{k,j,l}^{n-1} & \mbox{ if } & 1 \leq l \leq n-1 \\
r_k c_i^n & \mbox{ if } & l=n
\end{array} \right.$$

Note that $${\cal A}_k^{n} := \{ {\bf a_{k,j,i}^{n}} \}_{j,i} \subset S_{r_k}^{n-1} \subset \mathbb{R}^{n}$$  and therefore $${\cal A}^{n}:= \cup_{k=1}^p {\cal A}_k^{n}$$ has norm spectrum $R$.

Define also $$w^{n}:=w^{n-1} \star \gamma^n : {\cal A}^{n} \rightarrow \mathbb{R}^+;$$ recall that this means that $$w^{n}({\bf a_{k,j,i}^{n}})=w^{n}_{k,j,i}= \gamma_i^n w^{n-1}_{k,j}.$$  

Finally, note that, for every $1 \leq k \leq p$, we have $$W_k^{n}:=\sum_{j,i} w^{n}_{k,j,i}=\sum_{i} \gamma^{n}_{i} \cdot \sum_j w^{n-1}_{k,j}=\sum_j w^{n-1}_{k,j}=W_k^{n-1}.$$  

With these notations $({\cal A}^{n}, w^{n})$ is a Euclidean $t$-design, if for every $f \in \mathrm{Pol}_t(\mathbb{R}^{n})$, we have 

\begin{equation}\label{AA}
\sum_k \frac{W_k^{n}}{\sigma^{n}(S_{r_k}^{n-1})} \int_{S_{r_k}^{n-1}} f({\bf x}) d \sigma^{n}({\bf x}) = \sum_{k,j,i} w^{n}_{k,j,i} f({\bf a_{k,j,i}^{n}}).
\end{equation}

It suffices to verify (\ref{AA}) for monomials; for $$f({\bf x})=\prod_{l=1}^n x_l^{\alpha_l}$$ (\ref{AA}) becomes 

\begin{equation}\label{AB}
\sum_k \frac{W_k^{n}}{\sigma^{n}(S_{r_k}^{n-1})} \int_{S_{r_k}^{n-1}} \prod_{l=1}^n x_l^{\alpha_l} d \sigma^{n}({\bf x}) = \sum_{k,j,i} w^{n}_{k,j,i} \prod_{l=1}^n (a_{k,j,i,l}^{n})^{\alpha_l}.
\end{equation}

The key step in our proof is the identity

$$\int_{S^{n-1}} \prod_{l=1}^n x_l^{\alpha_l} d \sigma^{n}({\bf x}) = \int_{S^{n-2}} \prod_{l=1}^{n-1} x_l^{\alpha_l} d \sigma^{n-1}({\bf x}) \cdot \int_I (1-x^2)^{\frac{\sum_{l=1}^{n-1} \alpha_l}{2}} x^{\alpha_n} (1-x^2)^{\frac{n-3}{2}}  d x ;$$ see e.g. \cite{Str:1971a}.

By (\ref{GA}) and (\ref{GB}) we also have $$\sigma^{n}(S^{n-1})=\frac{2 \pi^{\frac{n}{2}}}{\Gamma(\frac{n}{2})}= \frac{2 \pi^{\frac{n-1}{2}}}{\Gamma(\frac{n-1}{2})} \cdot \frac{\sqrt{\pi} \cdot \Gamma(\frac{n-1}{2})}{\Gamma(\frac{n}{2})}= \sigma^{n-1}(S^{n-2}) \cdot \mu(I).$$

With these and other well known identities of spherical integration we re-write the left-hand side of (\ref{AB}), as follows.

\begin{eqnarray*}
\lefteqn{
\sum_k \frac{W_k^{n}}{\sigma^{n}(S_{r_k}^{n-1})} \int_{S_{r_k}^{n-1}} \prod_{l=1}^n x_l^{\alpha_l} d \sigma^{n}({\bf x})
} \\
& = & \sum_k \frac{W_k^{n}}{r_k^{n-1}\sigma^{n}(S_{1}^{n-1})} r_k^{n-1+\sum_{l=1}^n \alpha_l} 
\int_{S_{1}^{n-1}} \prod_{l=1}^n x_l^{\alpha_l}  d \sigma^{n}({\bf x}) \\
& = & \sum_k \frac{W_k^{n}}{r_k^{n-1}\sigma^{n-1}(S_{1}^{n-2}) \mu(I)} r_k^{n-1+\sum_{l=1}^n \alpha_l} \int_{S_{1}^{n-2}} \prod_{l=1}^{n-1} x_l^{\alpha_l}  d \sigma^{n-1}({\bf x}) \\ 
& & \cdot  \int_I (1-x^2)^{\frac{\sum_{l=1}^{n-1} \alpha_l}{2}} x^{\alpha_n} (1-x^2)^{\frac{n-3}{2}}  d x \\
& = & \sum_k \frac{W_k^{n}}{r_k\sigma^{n-1}(S_{r_k}^{n-2}) \mu(I)} r_k^{1+\alpha_n} \int_{S_{r_k}^{n-2}} \prod_{l=1}^{n-1} x_l^{\alpha_l}  d \sigma^{n-1}({\bf x}) \\ 
& & \cdot  \int_I (1-x^2)^{\frac{\sum_{l=1}^{n-1} \alpha_l}{2}} x^{\alpha_n} (1-x^2)^{\frac{n-3}{2}}  d x \\
& = & \sum_k \frac{W_k^{n}}{r_k\sigma^{n-1}(S_{r_k}^{n-2}) \mu(I)} r_k \int_{S_{r_k}^{n-2}} \left(\sum_{l=1}^{n-1} x_l^2 \right)^{\frac{\alpha_n}{2}}\prod_{l=1}^{n-1} x_l^{\alpha_l}  d \sigma^{n-1}({\bf x}) \\ 
& & \cdot  \int_I (1-x^2)^{\frac{\sum_{l=1}^{n-1} \alpha_l}{2}} x^{\alpha_n} (1-x^2)^{\frac{n-3}{2}}  d x \\
& = & \sum_k \frac{W_k^{n-1}}{\sigma^{n-1}(S_{r_k}^{n-2}) } \int_{S_{r_k}^{n-2}} \left(\sum_{l=1}^{n-1} x_l^2 \right)^{\frac{\alpha_n}{2}}\prod_{l=1}^{n-1} x_l^{\alpha_l}  d \sigma^{n-1}({\bf x}) \\ 
& & \cdot  \frac{1}{\mu(I)} \int_I (1-x^2)^{\frac{\sum_{l=1}^{n-1} \alpha_l}{2}} x^{\alpha_n} (1-x^2)^{\frac{n-3}{2}}  d x 
\end{eqnarray*}

We now re-write the right-hand side of (\ref{AB}):

\begin{eqnarray*}
\lefteqn{
\sum_{k,j,i} w^{n}_{k,j,i} \prod_{l=1}^n (a_{k,j,i,l}^{n})^{\alpha_l}
} \\
& = & \sum_{k,j,i} \gamma_i^n w^{n-1}_{k,j} \left(\prod_{l=1}^{n-1} \left(\sqrt{1-(c_i^n)^2} a_{k,j,l}^{n-1} \right)^{\alpha_l} \right) \cdot (r_k c_i^n)^{\alpha_n} \\
& = & \sum_{k,j} w^{n-1}_{k,j} \cdot (r_k)^{\alpha_n} \prod_{l=1}^{n-1} (a_{k,j,l}^{n-1})^{\alpha_l}  \\
& & \cdot \sum_{i} \gamma_i^n \left(1-(c_i^n)^2 \right)^{\frac{\sum_{l=1}^{n-1} \alpha_l}{2}}(c_i^n)^{\alpha_n} \\
& = & \sum_{k,j} w^{n-1}_{k,j} \cdot \left(\sum_{l=1}^{n-1} (a_{k,j,l}^{n-1})^2 \right)^{\frac{\alpha_n}{2}} \prod_{l=1}^{n-1} (a_{k,j,l}^{n-1})^{\alpha_l}  \\
& & \cdot \sum_{i} \gamma_i^n \left(1-(c_i^n)^2 \right)^{\frac{\sum_{l=1}^{n-1} \alpha_l}{2}}(c_i^n)^{\alpha_n}
\end{eqnarray*}

Therefore, (\ref{AB}) can be re-written in the form 
\begin{equation}\label{AC}
C(\mathbb{R}^{n-1}) \cdot C(I) = D(\mathbb{R}^{n-1}) \cdot D(I)
\end{equation}
where
\begin{eqnarray*}
C(\mathbb{R}^{n-1}) &= & \sum_k \frac{W_k^{n-1}}{\sigma^{n-1}(S_{r_k}^{n-2}) } \int_{S_{r_k}^{n-2}} \left(\sum_{l=1}^{n-1} x_l^2 \right)^{\frac{\alpha_n}{2}}\prod_{l=1}^{n-1} x_l^{\alpha_l}  d \sigma^{n-1}({\bf x}) \\
C(I) &= & \frac{1}{\mu(I)} \int_I (1-x^2)^{\frac{\sum_{l=1}^{n-1} \alpha_l}{2}} x^{\alpha_n} (1-x^2)^{\frac{n-3}{2}}  d x \\
D(\mathbb{R}^{n-1}) &= & \sum_{k,j} w^{n-1}_{k,j} \cdot \left(\sum_{l=1}^{n-1} (a_{k,j,l}^{n-1})^2 \right)^{\frac{\alpha_n}{2}} \prod_{l=1}^{n-1} (a_{k,j,l}^{n-1})^{\alpha_l} \\
D(I) &= & \sum_{i} \gamma_i^n \left(1-(c_i^n)^2 \right)^{\frac{\sum_{l=1}^{n-1} \alpha_l}{2}}(c_i^n)^{\alpha_n}
\end{eqnarray*}

We now distinguish three cases.

Case 1: $\alpha_n$ is odd.

Define the function $g: I \rightarrow \mathbb{R}$ by $$g(x)=(1-x^2)^{\frac{\sum_{l=1}^{n-1} \alpha_l}{2}} x^{\alpha_n} (1-x^2)^{\frac{n-3}{2}}.$$  Since $g$ is an odd function, $C(I)= 0$.  Also, since $({\cal C}^n, \gamma^n)$ is antipodal, $D(I) =0$.  Therefore (\ref{AC}) trivially holds as both sides equal zero.

Case 2: $\alpha_n$ is even and $\sum_{l=1}^{n-1} \alpha_l$ is odd.

Define the function $h: \mathbb{R}^{n-1} \rightarrow \mathbb{R}$ by $$h({\bf x}) = \left(\sum_{l=1}^{n-1} x_l^2 \right)^{\frac{\alpha_n}{2}}\prod_{l=1}^{n-1} x_l^{\alpha_l}.$$  Note that, since $\alpha_n$ is even, $h \in \mathrm{Pol}_t(\mathbb{R}^{n-1})$, and therefore we can apply (\ref{AD}) to get $C(\mathbb{R}^{n-1})=D(\mathbb{R}^{n-1})$.

But if $\sum_{l=1}^{n-1} \alpha_l$ is odd then at least one of these $n-1$ exponents is odd, so $C(\mathbb{R}^{n-1})=0$, and therefore $D(\mathbb{R}^{n-1})=0$.  This implies again that (\ref{AC}) trivially holds as both sides equal zero.

Case 3: $\alpha_n$ is even and $\sum_{l=1}^{n-1} \alpha_l$ is even.  

In this case, for the functions $g$ of Case 1 and $h$ of Case 2 we have $g \in \mathrm{Pol}_t(\mathbb{R})$ and $h \in \mathrm{Pol}_t(\mathbb{R}^{n-1})$, and therefore we can apply both (\ref{AE}) and (\ref{AD}).  We thus get $C(I)=D(I)$ and  $C(\mathbb{R}^{n-1})=D(\mathbb{R}^{n-1})$; multiplying these two equations yields (\ref{AC}), and this finishes the proof of Theorem \ref{FC}.

\vspace*{.2in}

\begin{center}

{\Large{{\bf Part II: Tight Euclidean Designs}}}

\end{center}

\section{Harmonic polynomials over $\mathbb{R}^n$}

There are several equivalent definitions of Euclidean designs.  For our purposes in this section and the next, the following will be convenient.

\begin{prop}[Neumaier--Seidel, \cite{NeuSei:1988a}]\label{FH}

The weighted set $({\cal X},w)$ is a Euclidean $t$-design in $\mathbb{R}^n$, if and only if,  
$$\sum_{{\bf x} \in {\cal X}} w({\bf x})||{\bf x}||^{2s_1} f({\bf x})=0$$ for every $0 \leq 2s_1 \leq t$ and $f \in  \mathrm{Harm}_s( \mathbb{R}^n)$ with $1 \leq s \leq t-2s_1$.

\end{prop}

In this section we review some information on harmonic polynomials over $\mathbb{R}^n$, and develop some useful results about a special subspace of $\mathrm{Harm}_s( \mathbb{R}^n)$.

Recall that a polynomial is \emph{harmonic} if it satisfies Laplace's equation $\Delta f=0$.  The set of homogeneous harmonic polynomials of degree $s$ in $\mathbb{R}^n$ forms the vector space $\mathrm{Harm}_{s}(\mathbb{R}^n)$ with $$\dim \mathrm{Harm}_{s}(\mathbb{R}^n)=\bc{n+s-1}{n-1}-\bc{n+s-3}{n-1}.$$

An explicit basis for $\mathrm{Harm}_{s}(\mathbb{R}^n)$ can be found, as follows.  

Choose integers $m_0, m_1, \dots, m_{n-2}$ with $0 \leq m_{n-2} \leq \cdots \leq m_1 \leq m_0=s$.  For $k=0,1,\dots, n-3$, set $$r_k=r_k(x_{k+1},\dots,x_{n}) =\sqrt{\sum_{l=k+1}^n x_l^2},$$ and define
$$g_k=g_k(x_{k+1},\dots,x_{n}) =r_k^{m_k-m_{k+1}} P_{m_k-m_{k+1}}^{m_{k+1}+(n-k-2)/2} \left(\frac{x_{k+1}}{r_k} \right)$$ where $P_{m_k-m_{k+1}}^{m_{k+1}+(n-k-2)/2}$ is the Gegenbauer polynomial introduced already in Section 2.  Note that $g_k$ is a polynomial of degree $m_k-m_{k+1}$.

Let also $$h_1(x_{n-1},x_n)=\mathrm{Re}(x_{n-1}+ix_n)^{m_{n-2}}$$ and $$h_2(x_{n-1},x_n)=\mathrm{Im}(x_{n-1}+ix_n)^{m_{n-2}}.$$

Finally, for integer(s) $1 \leq \mu \leq \mathrm{min} \{2, m_{n-2}+1\}$, define 
$$f_{m_0,m_1,\dots,m_{n-2},\mu}(x_1,\dots,x_n)=h_{\mu}(x_{n-1},x_n) \cdot \prod_{k=0}^{n-3}g_k(x_{k+1},\dots,x_{n})$$
and the set $$\Phi_s^n=\{f_{m_0,m_1,\dots,m_{n-2},\mu}(x_1,\dots,x_n) \}.$$
Note that $f_{m_0,m_1,\dots,m_{n-2},\mu}(x_1,\dots,x_n)$ is a polynomial of degree $s$ and that 
\begin{eqnarray*}
|\Phi_s^n| & = & |\{  (m_{n-2},\cdots, m_0,\mu) \in \mathbb{Z}^n \mbox{   } | \\ & & \mbox{             } \mbox{             } 0 \leq m_{n-2} \leq \cdots \leq m_0=s, 1 \leq \mu \leq \mathrm{min} \{2, m_{n-2}+1\} \}| \\
& = & 2 \cdot \bc{n+s-3}{n-2} + \bc{n+s-3}{n-3} \\
& = & \bc{n+s-1}{n-1}-\bc{n+s-3}{n-1} \\
& = & \dim \mathrm{Harm}_{s}(\mathbb{R}^n);
\end{eqnarray*} furthermore, we have the following.

\begin{prop}[\cite{Erd:1953a}]

With the above notations, the set 
$\Phi_s^n$
forms a basis for $\mathrm{Harm}_{s}(\mathbb{R}^n)$.
\end{prop}

In particular, we see that for $n=2$ we have $\dim \mathrm{Harm}_{s}(\mathbb{R}^2)=2$ and 
\begin{equation}\label{CC}
\Phi_s^2=\{h_1,h_2\};
\end{equation} 
we will use this in Section 5.  

For larger values of $n$, Proposition \ref{FH} is not convenient due to the large size of $\Phi_s^n$.  However, as we will soon see, if we consider only designs with a high degree of symmetry, the necessary criteria can be greatly reduced.  In particular, in Section 6, we will construct fully symmetric tight Euclidean designs.  A subset of $\mathbb{R}^n$ is said to be \emph{fully symmetric} if it remains fixed by (i) reflections with respect to co\"ordinate hyperplanes, and (ii) permutations of the co\"ordinates; a Euclidean design $({\cal X},w)$ is fully symmetric if ${\cal X}$ is a fully symmetric set, and $w(x)=w(y)$ whenever $x$ and $y$ are images of each other with respect to transformations of type (i) or (ii). 

Let us first address symmetry with respect to co\"ordinate hyperplanes.  For this latter purpose, we are specifically interested in \emph{fully even} harmonic polynomials, that is, those for which $$f(x_1, \dots, x_{l-1}, x_l, x_{l+1}, \dots, x_n)=f(x_1, \dots, x_{l-1}, -x_l, x_{l+1}, \dots, x_n)$$ holds for every $1 \leq l \leq n$; let $\mathrm{FEvenHarm}_{s}(\mathbb{R}^n)$ denote the set of fully even polynomials in $\mathrm{Harm}_{s}(\mathbb{R}^n)$.  We also let $\mathrm{POddHarm}_{s}(\mathbb{R}^n)$ consist of the \emph{partially odd} members of $\mathrm{Harm}_{s}(\mathbb{R}^n)$; these are polynomials with $$f(x_1, \dots, x_{l-1}, x_l, x_{l+1}, \dots, x_n)=-f(x_1, \dots, x_{l-1}, -x_l, x_{l+1}, \dots, x_n)$$ holds for some $1 \leq l \leq n$.  Note that in a polynomial which belongs to $\mathrm{FEvenHarm}_{s}(\mathbb{R}^n)$, in every term every variable has an even degree; while in members of $\mathrm{POddHarm}_{s}(\mathbb{R}^n)$, at least one variable appears only with an odd degree in every term.

We then let $$\mathrm{FEven}\Phi_s^n=\Phi_s^n \cap \mathrm{FEvenHarm}(\mathbb{R}^n),$$ and $$\mathrm{POdd}\Phi_s^n=\Phi_s^n \cap \mathrm{POddHarm}(\mathbb{R}^n).$$   

Recall that a Gegenbauer polynomial of degree $s$ is an even function when $s$ is even and an odd function when $s$ is odd.  Therefore the polynomial $g_k=g_k(x_{k+1},\dots,x_{n})$, defined above, has its variable $x_{k+1}$ with even exponents only if its degree $m_k-m_{k+1}$ is even, and odd exponents only if $m_k-m_{k+1}$ is odd; the variables $x_{k+2}, \dots, x_n$ all appear with even exponents only.  Consequently, we can determine easily which members of is in $\Phi_s^n$ are in $\mathrm{FEven}\Phi_s^n$, and we can count the elements of $\mathrm{FEven}\Phi_s^n$, as follows.

\begin{prop}\label{HA}  With the notations above, we have the following.
\begin{enumerate}
\item $$\mathrm\Phi_s^n = \mathrm{FEven}\Phi_s^n \cup \mathrm{POdd}\Phi_s^n;$$
\item A polynomial $f_{m_0,m_1,\dots,m_{n-2},\mu}(x_1,\dots,x_n) \in \Phi_s^n$ is in $\mathrm{E}\Phi_s^n$, if and only if, $m_{0},\cdots, m_{n-2}$ are all even and $\mu=1$;
\item
$\mathrm{FEvens}\Phi_s^n$ is a basis for $\mathrm{FEvenHarm}_{s}(\mathbb{R}^n)$;
\item If $s$ is even, then 
$$\dim \mathrm{FEvenHarm}_{s}(\mathbb{R}^n)=\bc{n+\frac{s}{2}-2}{n-2}.$$
\end{enumerate}
\end{prop}
In particular, for $s=2$ we have 
\begin{equation}\label{CA}
\dim \mathrm{FEvenHarm}_{2}(\mathbb{R}^n)=n-1,
\end{equation}
a substantial reduction compared to $\dim \mathrm{Harm}_{2}(\mathbb{R}^n)$.  In fact, it is easy to see that we can choose $$\mathrm{FEven}\Phi_2^n=\{f_2=x_i^2-x_{i+1}^2 \mbox{  } | \mbox{  } 1 \leq i \leq n-1 \}.$$  

For larger values of $s$, it proves even more profitable to additionally consider symmetry with respect to the permutation of co\"ordinates.  Consider first the combinatorial identity
\begin{equation}\label{CB}
\bc{n+\frac{s}{2}-2}{n-2}=\sum_{j=2}^{\frac{s}{2}}\bc{\frac{s}{2}-2}{j-2}\bc{n}{j},
\end{equation}
which holds for every even $s$ with $s \geq 4$.  

Based on (\ref{CB}), we will attempt to write $\mathrm{FEven}\Phi_s^n$ as the union of $$\sum_{j=2}^{\frac{s}{2}}\bc{\frac{s}{2}-2}{j-2}=2^{\frac{s}{2}-2}$$ families of functions, with each family of the form $$\{f_{s,j}(x_{i_1},\dots,x_{i_j}) \mbox{  } | \mbox{  } 1 \leq i_1 < \cdots < i_j \leq n \}$$ for some harmonic polynomial $f_{s,j}$ of degree $s$ and on $j$ variables with $2 \leq j \leq \frac{s}{2}$.  For our purposes in section 6, we will need to do this for $s=4$ and $s=6$; for these values we have 
\begin{eqnarray*} 
f_{4,2} & = & x_i^4-6x_i^2x_j^2+x_j^4, \\ 
f_{6,2} & = & x_i^6-15x_i^2x_j^2+15x_i^2x_j^2-x_j^6, \mbox{   and   } \\ 
f_{6,3} & = & 2(x_i^6+x_j^6+x_k^6)-15(x_i^4x_j^2+x_i^2x_j^4+x_i^4x_k^2+x_i^2x_k^4+x_j^4x_k^2+x_j^2x_k^4)+180x_i^2x_j^2x_k^2.
\end{eqnarray*}

Note further, that if $({\cal X},w)$ is fully symmetric, then the equation in Proposition \ref{FH} always holds for $f=f_2$ and $f=f_{6,2}$.  Therefore, we have the following.

\begin{cor}\label{FJ} Suppose that $({\cal X},w)$ is fully symmetric with norm spectrum $R$, and assume that its weight function is a constant $w_r$ on each layer ${\cal X}_r$ of ${\cal X}$.  Let $f_{4,2}$ and $f_{6,3}$ as above.  

\begin{enumerate}

\item $({\cal X},w)$ is at least a Euclidean $3$-design;

\item $({\cal X},w)$ is a Euclidean $5$-design, if and only if, 
\begin{equation}\label{DA}
\sum_{r \in R} \sum_{{\bf x} \in {\cal X}_r} w_r f_{4,2}({\bf x}) =0;
\end{equation}  

\item $({\cal X},w)$ is a Euclidean $7$-design, if and only if, 
\begin{equation}\label{DB}
 \left\{ 
\begin{array}{l}
\sum_{r \in R} \sum_{{\bf x} \in {\cal X}_r} w_r f_{4,2}({\bf x}) =0 \\ \\
\sum_{r \in R} \sum_{{\bf x} \in {\cal X}_r} w_r f_{6,3}({\bf x}) =0 \\ \\
\sum_{r \in R} \sum_{{\bf x} \in {\cal X}_r} w_r r^2 f_{4,2}({\bf x}) =0 
\end{array}
\right.
\end{equation}

\end{enumerate}

\end{cor}

\section{Tight Euclidean $t$-designs in the plane}

In this section we prove Theorem \ref{FD}.  We first recall the construction.

Let $1 \leq p \leq \left\lfloor \frac{t+5}{4}\right\rfloor$, and set $m=t+3-2p$.  Let $R=\{r_1,\dots,r_p\}$ be an arbitrary set of $p$ distinct positive real numbers; without loss of generality assume that $r_1<r_2< \cdots < r_p$.  For integers $k$ and $j$, define the point in the plane $$b_{k,j}= \left(r_k \cos \left( \frac{2j+k}{m} \pi \right), r_k \sin \left( \frac{2j+k}{m} \pi \right) \right),$$ and let $${\cal B} = \{ b_{k,j} \mbox{  } | \mbox{  } 1 \leq j \leq m , 1 \leq k \leq p \}.$$   Then ${\cal B} $ is the union of $p$ concentric regular $m$-gons centered at the origin.  Note that when $t$ is odd then $m$ is even, and therefore the set ${\cal B}$ is antipodal (but ${\cal B} \cap -{\cal B} = \emptyset$ when $t$ is even).

Next, define the weight function $w : {\cal B} \rightarrow \mathbb{R}^+$ with
$$w(b_{k,j})=\left\{ \begin{array}{lll}
\frac{1}{r_1^m} & \mbox{ if } & k=1 \\
(-1)^k \frac{1}{r_k^m} \prod_{2 \leq l \leq p, l \not = k} \frac{r_1^2-r_l^2}{ r_k^2-r_l^2} & \mbox{ if } & 2 \leq k \leq p
\end{array} \right.$$

Note that $w(b_{k,j})$ is always positive.

In order to prove that $({\cal B},w)$ is a tight Euclidean $t$-design, we first verify that $N(2,p,t)=|{\cal B}|$.  Clearly, we have  $|{\cal B}|=p \cdot m=p(t+3-2p)$.  On the other hand, we have 
\begin{eqnarray*}
N(2,p,t) &= & \sum_{k=1}^p N_k \\
& = & \sum_{k=1}^p \bc{\left\lfloor \frac{t}{2} \right\rfloor +2-2k+1}{1}+\bc{\left\lfloor \frac{t-1}{2} \right\rfloor +2-2k+1}{1} \\
& = & p(t+3-2p),
\end{eqnarray*}
as claimed.

We need to prove that $({\cal B},w)$ is a Euclidean $t$-design.  As we saw in the last section, for $n=2$ we have $\dim \mathrm{Harm}_s(\mathbb{R}^2)=2$; furthermore, we found in (\ref{CC}) that the polynomials $\mathrm{Re}(x+iy)^s$ and $\mathrm{Im}(x+iy)^s$ form a basis for $\mathrm{Harm}_s(\mathbb{R}^2)$ (we will not distinguish between the complex number $x+iy$ and the point $(x,y)$).  Therefore, we see that $({\cal X},w)$ is a Euclidean $t$-design in $\mathbb{R}^2$, if and only if, $$\sum_{z \in {\cal X}} w(z) ||z||^{2s_1} z^s=0$$ for every $1 \leq s \leq t$ and $0 \leq s_1 \leq \lfloor \frac{t-s}{2} \rfloor$.

In view of this, we turn to complex numbers.
Let $$\psi = e^{\frac{\pi i}{m}} = \cos (\frac{\pi}{m} ) + i \sin(\frac{\pi}{m})$$ and $\lambda=\psi^2$.  Then $b_{k,j}=r_k \cdot \psi^{k} \cdot \lambda^j$.

To prove that $({\cal B}, w)$ is a Euclidean $t$-design, we need to show that for every $1 \leq s \leq t$ and $0 \leq s_1 \leq \lfloor \frac{t-s}{2} \rfloor$, we have $$\sum_{k=1}^p w_k \cdot r_k^{2s_1+s} \cdot \psi^{ks} \cdot \sum_{j=1}^m \lambda^{js} = 0.$$ 

Now $\lambda^s=1$ if and only if $m$ is an integer divisor of $s$, and so $$\sum_{j=1}^m \lambda^{js} =\left\{ \begin{array}{llll}
\lambda^{s} \cdot \frac{\lambda^{sm}-1}{\lambda^{s}-1} & = 0 & \mbox{ if } & m \not | s \\
\sum_{j=1}^m 1 & = m & \mbox{ if } & m | s
\end{array} \right.$$

Furthermore, $$m=t+3-2p \geq t+3 -2 \cdot \frac{t+5}{4} = \frac{t+1}{2} > \frac{s}{2},$$ so $m$ can only be a divisor of $s$ when $m=s$.

Therefore, we only need to show that 
\begin{equation}\label{BB}
\sum_k w_k (-1)^k r_k^{2s_1+m}=0
\end{equation} 
holds for every $0 \leq s_1 \leq \lfloor (t-m)/2 \rfloor =p-2$ (note that $\psi^{km}=(-1)^k$).

We can re-write (\ref{BB}) as 
$$\sum_{k=2}^p w_k (-1)^k r_k^{2s_1+m}=w_1 r_1^{2s_1+m},$$
which, after substituting the value of the weights, becomes

$$\sum_{k=2}^p r_k^{2s_1} \frac{\prod_{2 \leq l \leq p, l \not = k} (r_1^2-r_l^2)}{\prod_{2 \leq l \leq p, l \not = k} (r_k^2-r_l^2)}= r_1^{2s_1},$$
which indeed holds for every $s_1=0,1,2, \dots, p-2$ according to Vandermonde's Theorem.  This completes our proof.

\section{Tight designs in dimension $n \geq 3$}

In this section we construct examples for tight Euclidean $t$-designs for $t=5$ and $t=7$. 

We look for Euclidean designs inside the integer lattice $\mathbb{Z}^n$; in particular, here we restrict our attention to points within the box $${\cal I}^n = I^n \cap \mathbb{Z}^n = \{ {\bf x}=(x_1,\dots,x_n) \mbox{  } | \mbox{  } x_i \in \{-1,0,1\} \}.$$  We partition the $3^n$ elements in ${\cal I}^n$ according to their norm.  Namely, we let $${\cal I}^n_k =\{ {\bf x} \in {\cal I}^n \mbox{   } | \mbox{  } ||{\bf x} ||^2=k \}.$$  Note that $$|{\cal I}^n_k|=2^k \cdot \bc{n}{k}.$$  For example, for $n=3$, this partition is the following:  $${\cal I}^3 =\{ {\bf 0} \} \cup \{ \mbox{octahedron} \}  \cup \{ \mbox{cuboctahedron} \} \cup \{ \mbox{cube} \}.$$

We consider sets ${\cal X}$ of the form $${\cal X}={\cal X}(J)=\cup_{k \in J} \frac{r_k}{\sqrt{k}} {\cal I}^n_k$$ where $J \subset \{1,2,\dots,n\}$ and $r_k>0$ for every $k \in J$.  Note that ${\cal X}$ is antipodal and has norm spectrum $R=\{r_k \mbox{  } | \mbox{  } k \in J \}$.  The weight function $w: {\cal X} \rightarrow \mathbb{R}^+$ on ${\cal X}$ will be constant on each layer of ${\cal X}$ (see \cite{BanBan:Prepa}); let us denote the weight of ${\bf x} \in \frac{r_k}{\sqrt{k}}{\cal I}^n_k$ by $w({\bf x})=w_k$.  

Our goal is to find index sets $J$ and positive numbers $r_k$ and $w_k$ ($k \in J$) for which $({\cal X},w)$ is a tight Euclidean design.   Clearly, without loss of generality (see \cite{BanBan:Prepa}) we can choose one of the radii and one of the weights freely; since in our examples we will always have $1 \in J$, we let $r_1=1$ and $w_1=1$.

Note that $({\cal X},w)$ is fully symmetric, that is, it remains fixed by permutations of the co\"ordinates and reflections with respect to co\"ordinate hyperplanes.  Therefore, to ascertain that it is a Euclidean $t$-design, it suffices to use the techniques of section 4; in particular, in the case of $t \leq 7$, we may use Corollary \ref{FJ}.

For a given $f \in \Phi_s^n$, let $$S(f)=\sum_{{\bf x} \in {\cal I}^n_k} f({\bf x}).$$  We can compute that
$$S(f_{4,2})= 2\cdot 2^k \cdot \bc{n-1}{k-1} - 6 \cdot 2^k \cdot \bc{n-2}{k-2}$$ and
$$S(f_{6,3})= 6\cdot 2^k \cdot \bc{n-1}{k-1} - 90 \cdot 2^k \cdot \bc{n-2}{k-2}+180 \cdot 2^k \cdot \bc{n-3}{k-3}.$$  Now we are ready to discuss some specific examples and to prove Proposition \ref{FG}.

\subsection{$n=3$, $t=5$}

Note that for $p \geq 2$, we have
$N(3,p,5)=14$.  Since $|{\cal I}^3_1|=6$ and $|{\cal I}^3_3|=8$, we choose $J=\{1,3\}$.  We then have $${\cal X}={\cal I}^3_1 \cup \frac{r_3}{\sqrt{3}} {\cal I}^3_3.$$    

Equation (\ref{DA}) then becomes
$$4-\frac{32}{9} r_3^4w_3=0,$$
and this holds exactly when $$w_3=\frac{9}{8r_3^4},$$ proving Proposition \ref{FG} (1).

\subsection{$n=3$, $t=7$}

Note that for $p \geq 2$, we have
$N(3,2,7)=26$.  Since $|{\cal I}^3_1|=6$,  $|{\cal I}^3_2|=12$, and $|{\cal I}^3_3|=8$, we choose $J=\{1,2,3\}$.  We then have $${\cal X}={\cal I}^3_1 \cup \frac{r_2}{\sqrt{2}} {\cal I}^3_2\cup \frac{r_3}{\sqrt{3}} {\cal I}^3_3.$$    

The system of equations (\ref{DB}) becomes

$$\left\{ 
\begin{array}{l}
4-2 r_2^4w_2-\frac{32}{9} r_3^4w_3 =0 \\ \\
 12-39r_2^6w_2+\frac{256}{9} r_3^6w_3 =0 \\ \\
4-2 r_2^6w_2-\frac{32}{9} r_3^6w_3 =0 
\end{array}
\right.$$

All three equations hold, if and only if, $$r_2=\sqrt{\frac{2}{5r_3^2-3}}r_3, \mbox{   } w_2=\frac{(5r_3^2-3)^3}{10r_3^6}, \mbox{   and   } w_3=\frac{27}{40r_3^6}.$$

Letting $\lambda=\frac{5r_3^2-3}{2}$ produces Proposition \ref{FG} (2).

Note that $|\{1,r_2,r_3\}| \leq 2$ implies that $1=r_2=r_3$, contradicting $p \geq 2$; therefore we must have $p=3$.  Two particularly nice choices for the parameters are worth mentioning:
$${\cal X}={\cal I}^3_1 \cup {\cal I}^3_2\cup \frac{1}{2} {\cal I}^3_3$$ with weights $w_2=\frac{1}{10}$ and $w_3=\frac{8}{5}$; and
$${\cal X}={\cal I}^3_1 \cup \frac{1}{2}{\cal I}^3_2\cup {\cal I}^3_3$$ with weights $w_2=\frac{32}{5}$ and $w_3=\frac{1}{40}$.

\subsection{$n=4$, $t=7$}

Note that for $p \geq 2$, we have
$N(4,2,7)=48$.  Since $|{\cal I}^4_1|=8$,  $|{\cal I}^4_2|=24$, and $|{\cal I}^4_4|=16$, we choose $J=\{1,2,4\}$.  We then have $${\cal X}={\cal I}^4_1 \cup \frac{r_2}{\sqrt{2}} {\cal I}^4_2\cup \frac{r_4}{2} {\cal I}^4_4.$$    

The system of equations (\ref{DB}) becomes

$$\left\{ 
\begin{array}{l}
4-4 r_4^4w_4 =0 \\ \\
 12-36r_2^6w_2+24r_4^6w_4 =0 \\ \\
4-4r_4^6w_4 =0 
\end{array}
\right.$$

All three equations hold, if and only if, $$r_4=1, \mbox{   } w_2=\frac{1}{r_2^6}, \mbox{   and   } w_4=1.$$

Therefore, we have $${\cal X}=({\cal I}^4_1 \cup \frac{1}{2} {\cal I}^4_4) \cup \frac{r_2}{\sqrt{2}} {\cal I}^4_2,$$
where ${\cal I}^4_2$ is the set of vectors with minimal (non-zero) norm in the ``checkerboard'' lattice $D_4$, and ${\cal I}^4_1 \cup \frac{1}{2} {\cal I}^4_4$ is the set of vectors with minimal (non-zero) norm in its dual lattice $D_4^*$.

This completes the proof of Proposition \ref{FG}.

{\bf Acknowledgments.}  Work on this paper began during a workshop at the Universit\'e de Gen\`eve, organized by Pierre de la Harpe during May 2004; I am very grateful for his hospitality and continued support.  I would also like to thank Eiichi Bannai for many helpful discussions.


\begin{thebibliography}{10}



\bibitem{Baj:1991d}
B.~Bajnok.
\newblock {C}hebyshev-type quadrature formulas on the sphere.
\newblock {\em Congr. Numer.}, 85:214--218, 1991.

\bibitem{Baj:1992a}
B.~Bajnok.
\newblock Construction of spherical $t$-designs.
\newblock {\em Geom. Dedicata}, 43/2:167--179, 1992.

\bibitem{Ban:1988a}
E.~Bannai.
\newblock On extremal finite sets in the sphere and other metric spaces.
\newblock {\em London Math. Soc. Lecture Note Ser.}, 131:13--38, 1988.

\bibitem{BanBan:Prepa}
E.~Bannai and E.~Bannai.
\newblock On Euclidean tight 4-designs.
\newblock Preprint.

\bibitem{BanBanSta:1983a}
E.~Bannai, E. Bannai,, and D. Stanton.
\newblock An upper bound for the cardinality of an $s$-distance subset in real Euclidean space {II}.
\newblock {\em Combinatorica}, 3/2:147--152, 1983.

\bibitem{BanDam:1979a}
E.~Bannai and R.~M. Damerell.
\newblock Tight spherical designs {I}.
\newblock {\em J. Math. Soc. Japan}, 31/1:199--207, 1979.

\bibitem{BanDam:1980a}
E.~Bannai and R.~M. Damerell.
\newblock Tight spherical designs {II}.
\newblock {\em J. London Math. Soc. (2)}, 21/1:13--30, 1980.

\bibitem{BanMunVen:Prepa}
E.~Bannai, A. Munemasa, and B. Venkov.
\newblock The nonexistence of certain tight spherical designs.
\newblock Preprint.

\bibitem{ConSlo:1999a}
J. H. Conway and N. J. A. Sloane.
\newblock {\em Sphere Packings, Lattices and Groups, 3rd ed}.
\newblock Springer--Verlag, 1999.

\bibitem{DelPac:Prepa}
P. de la Harpe and C. Pache.
\newblock Cubature formulas, geometrical designs, reproducing kernels, and Markov operators.
\newblock Preprint.

\bibitem{DelGoeSei:1977a}
P.~Delsarte, J.~M. Goethals, and J.~J. Seidel.
\newblock Spherical codes and designs.
\newblock {\em Geom. Dedicata}, 6/3:363--388, 1977.

\bibitem{DelSei:1989a}
P.~Delsarte and J.~J. Seidel.
\newblock Fisher type inequalities for Euclidean $t$-designs.
\newblock {\em Linear Algebra Appl.}, 114--115:213--230, 1989.

\bibitem{Erd:1953a}
A. Erd\'elyi et al.
\newblock Higher transcendental functions, Vol. II. (Bateman Manuscript Project).
\newblock MacGraw--Hill, 1953.

\bibitem{EriZin:2001a}
T. Ericson and V. Zinoviev.
\newblock Codes on Euclidean spheres.
\newblock Elsevier, 2001.

\bibitem{God:1993a}
C.~D. Godsil.
\newblock {\em Algebraic Combinatorics}.
\newblock Chapman and Hall, Inc., 1993.

\bibitem{GoeSei:1979a}
J.~M. Goethals and J.~J. Seidel.
\newblock Spherical designs.
\newblock In D.~K. Ray-Chaudhuri, editor, {\em Relations between combinatorics
  and other parts of mathematics}, 
\newblock {\em Proc. Sympos. Pure Math.}, 34:255--272. Amer. {M}ath. {S}oc., {P}rovidence, {RI}, 1979.

\bibitem{GoeSei:1981a}
J.~M. Goethals and J.~J. Seidel.
\newblock Cubature formulae, polytopes and spherical designs.
\newblock In C.~Davis, B.~Gr{\"{u}}nbaum, and F.~A. Sher, editors, {\em The
  Geometric Vein: The Coxeter Festschrift}, pp. 203--218. Springer-{V}erlag
  New York -- Berlin, 1981.

\bibitem{Hil:1974a}
F. B. Hildebrand.
\newblock {\em Introduction to Numerical Analysis, 2nd ed}.
\newblock McGraw--Hill, Inc., 1974.

\bibitem{KorMey:1994a}
J.~Korevaar and J.~L.~H. Meyers.
\newblock Chebyshev-type quadrature on multidimensional domains.
\newblock {\em J. Approx. Theory}, 79:144--164, 1994.

\bibitem{Kui:1995a}
A. Kuijlaars.
\newblock Chebyshev-type quadrature for Jacobi weight functions.
\newblock {\em J. Comput. Applied Math.}, 57:171--180, 1995.

\bibitem{Kry:1962a}
V. I. Krylov.
\newblock {\em Approximate Calculation of Integrals}.
\newblock The Macmillan Company, New York, 1962.

\bibitem{NeuSei:1988a}
A. Neumaier and J. J. Seidel.
\newblock Discrete measures for spherical designs, eutactic stars and lattices.
\newblock {\em Nederl. Akad. Wetensch. Indag. Math.}, 50/3:321-334, 1988.

\bibitem{Sei:1990a}
J.~J. Seidel.
\newblock Designs and approximation.
\newblock {\em Contemp. Math.}, 111:179--186, 1990.

\bibitem{Sei:1994a}
J.~J. Seidel.
\newblock Isometric embeddings and geometric designs.
\newblock {\em Discrete Math.}, 136:281--293, 1994.

\bibitem{Sei:1996a}
J.~J. Seidel.
\newblock Spherical designs and tensors.
\newblock In E.~Bannai and A.~Munemasa, editors, {\em Progress in algebraic
  combinatorics}, 
\newblock {\em Adv. Stud. Pure Math.}, 24:309--321. Math. {S}oc. {J}apan, {T}okyo, 1996.

\bibitem{SeyZas:1984a}
P.~D. Seymour and T.~Zaslavsky.
\newblock Averaging sets: A generalization of mean values and spherical designs.
\newblock {\em Adv. Math.}, 52:213--240, 1984.

\bibitem{Str:1971a}
A. H. Stroud.
\newblock {\em Approximate Calculation of Multiple Integrals}.
\newblock Prentice--Hall, Inc., 1971.


\end{thebibliography}
\end{document}